\documentclass{amsart}
\usepackage{amssymb,graphicx,url}
\usepackage[line,ps,poly,all]{xy}
\newcommand{\ZZ}{\mathbb Z}%
\newcommand{\QQ}{\mathbb Q}%
\newcommand{\CC}{\mathbb C}%
\newcommand{\RR}{\mathbb R}%
\newcommand{\OO}{{\mathcal O}}

\newcommand{\G}{\Gamma}%
\newcommand{\Gb}{\bar{\Gamma}}%
\newcommand{\GC}{{\mathcal G}}%
\newcommand{\FF}{{\mathcal F}}%
\newcommand{\BC}{\mathcal{B}}%
\newcommand{\p}{\mathfrak{p}}%
\newcommand{\lsp}[1]{{}^{#1}\!}%

\DeclareMathOperator{\Stab}{Stab}%

\newcommand{\Par} {{\mathcal P}}%
\newcommand{\PP} {\mathbb{P}}%

\DeclareMathOperator{\SU}{SU}
\DeclareMathOperator{\SL}{SL}%

\newcommand{\I} {{\mathcal I}}%

\newcommand{\QQC}{\mathcal Q}%
\newcommand{\J} {{\mathcal J}}%
\newcommand{\D} {{\mathcal D}}%
\renewcommand{\Re}{\operatorname{Re}}%
\renewcommand{\Im}{\operatorname{Im}}%

\newcommand{\vect}[1]{\begin{pmatrix} #1  \end{pmatrix}}%
\newcommand{\mat}[1]{\begin{pmatrix} #1 \end{pmatrix}}%

\theoremstyle{plain}
\newtheorem{thm}{Theorem}[section]

\newtheorem{prop}[thm]{Proposition}

\theoremstyle{definition}
\newtheorem{defn}[thm]{Definition}

\newcommand{\group}[1]{\mathbf{#1}}

\begin{document}


\title{Integral cohomology of certain Picard modular surfaces}
\author{Dan Yasaki}
\address{Department of Mathematics and Statistics\\Lederle Graduate Research Tower\\ University of Massachusetts\\Amherst, MA 01003-9305}
\email{yasaki@math.umass.edu}
\date{}
\thanks{The original manuscript was prepared with the \AmS-\LaTeX\ macro
system and the \Xy-pic\ package.}
\keywords{spine, Picard modular group, locally symmetric space, cohomology of arithmetic subgroups}
\subjclass[2000]{Primary 11F75}
\begin{abstract}
Let $\Gb$ be the Picard modular group of an imaginary quadratic number field $k$ and let $\D$ be the associated symmetric space.  Let $\G \subset \Gb$ be a congruence subgroup.  We describe a method to compute the integral cohomology of the locally symmetric space $\G \backslash \D$.  The method is implemented for the cases $k=\QQ(i)$ and $k=\QQ(\sqrt{-3})$, and the cohomology is computed for various $\G$.  
\end{abstract}

\maketitle
\bibliographystyle{../../amsplain_initials.bst}

\begin{section}{Introduction}\label{sec:introduction}
Let $\Gb=\SU(2,1;\OO_k)$ be the Picard modular group of an imaginary quadratic number field $k$ and let $\D$ be the associated symmetric space.  Let $\G \subset \Gb$ be a congruence subgroup.  Although $\D$ is $4$-dimensional, the virtual cohomological dimension of $\Gb$ is $3$.  Hence the cohomology $H^i(\G \backslash \D)$ vanishes for $i > 3$.  

If $\G$ is torsion-free, then $\G \backslash \D$ is an Eilenberg-MacLane space for $\G$.  It follows that the group cohomology of $\G$ with trivial complex coefficients is isomorphic to the complex cohomology of the locally symmetric space,
\begin{equation}\label{eq:coh}
H^*(\G;\CC) \simeq H^*(\G \backslash \D;\CC).
\end{equation}
In fact, \eqref{eq:coh} remains true while using complex coefficients when $\G$ has torsion, but is not true in general for integral coefficients.  

A deep result of Franke \cite{Fra} describes a relationship between cohomology and automorphic forms.  The cohomology groups $H^*(\G;\CC)$, or more generally $H^*(\G;M)$ for any complex finite dimensional rational representation of $\SU(2,1)$, decompose into cuspidal cohomology and Eisenstein cohomology.  The cuspidal cohomology of $\G$ can be represented by cuspidal automorphic forms.  It is possible to compute the space of cuspidal Picard modular forms with the Hecke action using the Jacquet-Langlands correspondence \cite{De}.  However, these methods will not compute the torsion classes in the integral cohomology of the group or the locally symmetric space.  In this paper, we use topological methods to compute the torsion classes in the integral cohomology of the locally symmetric space $\G \backslash \D$.
 
There is a $3$-dimensional cell complex $W$, known as a spine, that can be used to compute the integral cohomology of the locally symmetric space $\G \backslash \D$.  The existence of such a $W$ for general $\QQ$-rank $1$ groups is known \cite{Yasrank1}, but there are few explicit examples for non-linear symmetric spaces.  We outline a  method of computing $W$ for $\SU(2,1;\OO_k)$, where $k$ is an imaginary quadratic number field with class number $1$.  The structure of the spine is computed in \cite{Yaspicard} for the Gaussian case $k=\QQ(i)$ and is computed here for the Eisenstein case $k=\QQ(\sqrt{-3})$.  Falbel and Parker \cite{FaPa} do similar computations for the Eisenstein-Picard modular group, but with a different purpose.  They exhibit fundamental domain for the action of $\Gb$ on $\D$.  Using the the structure and combinatorics of the fundamental domain, they deduce a presentation for $\Gb$.    

The outline of the paper is as follows.  We first recall the Picard modular group and associated symmetric space in Section~\ref{sec:prelim}.  A method of computing  the spine $W$ is given in Section~\ref{sec:spine}.  Section~\ref{sec:cohomology} outlines the method of \cite{AGM} and its implications in the context of our cell complex.  The cell complex is computed in \cite{Yaspicard} for the Picard modular group over the Gaussian integers. We compute the cell complex and give the stabilizers for Picard modular group over the Eisenstein integers in Section~\ref{sec:examples}.  Finally, Section~\ref{sec:numerical} gives the cohomology computation results.

I would like to thank Paul Gunnells for many helpful comments and patiently explaining his paper to me.  I also thank T.N. Venkataramana, who helped with understanding \cite{BlRo}. 
\end{section}

\begin{section}{Preliminaries}\label{sec:prelim}
Let $k$ be an imaginary quadratic field with discriminant $D$ and ring of integers $\OO$.  Thus $k=\QQ(\sqrt{D})$ and $D$ is either square-free and $D\equiv 1 \bmod 4 $ or $D=4D'$, where $D'$ is square-free and $D' \equiv 2 \bmod 4 $ or $D' \equiv 3 \bmod{4}$.  Then $\OO$ is generated by $1$ and $\omega = (D+\sqrt{D})/2$.  For $d=1,2,3,7,11,19,43,67,163$,  $\QQ(\sqrt{-d})$ has class number $h(k)= 1$ and $\OO$ is a principal ideal domain.  Fix and imaginary quadratic field $k$ with class number $1$.

\begin{subsection}{The unitary group}
Let $V$ be a $3$-dimensional $k$-vector space with an integral structure given by an $\OO$-lattice $L \subset V$.  Let ${\mathcal Q}:V \times V \to k$ be a non-degenerate Hermitian form on $V$ which is $\OO$-valued on $L$ and whose signature is $(2,1)$ on $V_\RR$.  Then $\G=\SU({\mathcal Q},V)$ is a semisimple algebraic group defined over $\QQ$ whose group of real points is isomorphic to $\SU(2,1)$.  The \emph{Picard Modular Group} of $k$ is 
\[\Gb=\group{G}(\ZZ)=\{\gamma \in \group{G}(\QQ)\; |\; \gamma L=L\}.\]

For this paper, fix an embedding of $k$ in $\CC$ and a $k$-basis of $V$ with respect to which $\mathcal{Q}$ is represented by
\begin{equation*}
\mathcal{Q}(u,v)=u^*C v, \quad \text{where $C=\begin{pmatrix} 0&0&\sqrt{D}^{-1}\\0&1&0\\-\sqrt{D}^{-1}&0&0 \end{pmatrix}$}.
\end{equation*}  
In particular,
\[
G=\group{G}(\RR)=\left\{g \in \SL(3,\CC)\; |\; g^*C g=C \right\}.\] 

\end{subsection}

\begin{subsection}{Symmetric space}
Let $\theta$ denote the Cartan involution given by inverse conjugate transpose and let $K$ be the fixed points under $\theta$.   Let $\D=G/K$ be the associated Riemannian symmetric space of non-compact type.  The symmetric space $\D=G/K$ has many useful realizations.

Using horospherical coordinates $(y,\beta,r)$ we can view the symmetric space as $\RR_{>0} \times \CC \times \RR$ as follows.  Let $P_0 \subseteq G$ be the rational parabolic subgroup of upper triangular matrices.
\begin{align*}
P_0&=\left\{\left.\begin{pmatrix} y\zeta & \beta \zeta^{-2} & \frac{\zeta}{y}\left(r - \frac{|\beta|^2}{2\sqrt{D}}\right) \\ 0 & \zeta^{-2} & -\frac{{\overline \beta}\zeta}{ \sqrt{D} y} \\ 0 & 0 & \zeta/y \end{pmatrix} \; \right| \;  \zeta,\beta \in \CC,\ |\zeta|=1,\ r \in \RR,\ y\in \RR_{>0} \right\}.
\end{align*}
$P_0$ acts transitively on $\D$, and every point $X \in \D$ can be written as $p  K$ for some $p \in P_0$.  When $p$ is written as above, the point $X=pK$ is independent of $\zeta$, and so we will denote such a point $X=(y,\beta,r)$.  
  
A line $l$ in $V$ is said to be \emph{negative} if 
\[\QQC(v,v) < 0, \quad \text{where $v$ is any vector on $l$.}\]
Let $\mathcal{L}$ denote the set negative lines in $V$.  The group $G$ acts transitively on $\mathcal{L}$, and the stabilizer of a negative line is a maximal compact subgroup of $G$.  Thus we can identify $\D$ with $\mathcal{L}$.  A negative line is the span of a vector with non-zero third component, and so each line can be identified with a vector $(z,u,1)^t$.  This identification gives the \emph{Siegel domain} realization of $\D$,
\[\D=\left\{(z,u) \in \CC^2\; \left|\; |u|^2 < \frac{2\Im(z)}{\sqrt{|D|}}\right.\right\}.\]
In coordinates, one computes that
\begin{equation}\label{eq:coordinates}
z=r+iy^2 -\frac{|\beta|^2}{2\sqrt{D}}\quad \text{and} \quad u=\frac{-\bar{\beta}}{\sqrt{D}}.
\end{equation}
\end{subsection}

\begin{subsection}{The cusps}
The cusps correspond to rational isotropic lines in $V$.  These lines are in 1-1 correspondence with the rational parabolic subgroups.  In particular, every rational parabolic subgroup stabilizes a rational isotropic line.  We associate to each rational parabolic subgroup $P$ an isotropic vector $v_P=(n,p,q)^t \in \OO^3$ from the line stabilized by $P$ such that the ideal $(n,p,q)=\OO$.  Note that this vector is well-defined up to multiplication $\OO^*$.  Thus for the rest of the paper, we identify vectors in $\OO^3$ that differ by $\OO^*$.  The isotropic condition $ \mathcal{Q}(v_P,v_P)=0$ implies that $v_P=(n,p,q)^t$ satisfies 
\begin{equation}\label{eq:isotropic}
|p|^2=\frac{1}{\sqrt{D}}(n \overline{q}-\overline{n}q).
\end{equation}  

Let $\Par$ denote the set of proper rational parabolic subgroups of $G$.  These parabolic subgroups are conjugate to $P_0$ via elements of $\group{G}(\QQ)$.  There is a natural action of $\Gb$ on $\Par$ given by conjugation and denoted
\[
\lsp{g}P :=gPg^{-1}. 
\]
The quotient $\Gb \backslash \Par$ is a finite set whose cardinality is the \emph{class number of $\Gb$}.  Zink has shown that $h(\Gb)=h(k)$ \cite{Z}.  In particular, since $h(k)=1$, every rational parabolic subgroup is conjugate to $P_0$ via an element of $\Gb$. 
\end{subsection}

\begin{subsection}{The spine}
There is a $\Gb$-invariant decomposition of $\D$ into codimension $0$ sets using exhaustion functions.  The union $\D_0$ of the boundaries of these sets forms a contractible $\Gb$-equivariant retraction of $\D$ that is known as a \emph{spine}.  This construction is described for the general $\QQ$-rank 1 case in \cite{Yasrank1} and is analogous to Ash's \emph{well-rounded retraction} \cite{A2} for linear symmetric spaces.

There is a collection of subsets of rational parabolic subgroups $\mathcal S$ called \emph{strongly admissible subsets} that give a decomposition of the spine 
\[\D_0=\coprod_{\substack{\I \in {\mathcal S}\\ |\I| >1}} \D'(\I).\] 
The group $\Gb$ preserves this decomposition.  In particular, 
\[\gamma \cdot \D'(\I) = \D'(\lsp{\gamma}\I)\]
for every $\gamma \in \Gb$.

The decomposition of $\D_0$ may be refined to a regular cell complex $W$ that has the property that the stabilizer of each cell fixes the cell pointwise.  Thus, the cohomology of $\Gb\backslash \D$ can be described by finite combinatorial data.
\end{subsection}

\begin{subsection}{Congruence subgroups}
Let $N\in \OO$ and let $\bar{k}=\OO/N\OO$.  Let $\bar{V}=\bar{k}^3$ and $\PP^2=\PP(\bar{V})$.  By $\PP^2$, we mean the set of vectors $(x_1,x_2,x_3) \in V$ that are primitive in the sense that the ideal $(x_1,x_2,x_3) =\OO$.  When $N$ is prime, this is the usual projective space over the residue field $\bar{k}$.  The equivalence class of the vector $(x_1,x_2,x_3)$ will be denoted $[x_1:x_2:x_3]$.  We view these triples as row vectors, and so $\Gb$ acts on $\bar{V}$ and $\PP^2$ on the right.

Let $\G_1(N)\subset \Gb$ denote the stabilizer of the vector $(0,0,1) \in V$, and let $\G_0(N)\subset \Gb$ denote the stabilizer of the point $[0:0:1]\in \PP$.  Let $\G(N)$ denote the kernel of the map $\mu:\Gb \to \SL_3(\bar{k})$.  
\end{subsection}
\end{section}

\begin{section}{Computing the spine}\label{sec:spine}
In any given example, one must understand the exhaustion functions to compute the spine.  The functions are parametrized by rational parabolic subgroups can be thought of as height functions with respect to the various cusps.  The family of exhaustion functions are $\Gb$-invariant in the sense that 
\begin{equation}\label{eq:invariantfamily}
f_{\lsp{\gamma}P}(X)=f_P(\gamma^{-1}\cdot X)\quad \text{for $\gamma \in \Gb$.}
\end{equation}
These exhaustion functions for $k=\QQ(i)$ are described in detail in \cite{Yaspicard}.  Here, we give the exhaustion functions in coordinates for an imaginary quadratic field $k=\QQ(\sqrt{D})$ with class number one.

\begin{subsection}{The exhaustion functions}
Let $X=(z,u)\in \D$ and let $P$ be a rational parabolic subgroup of $G$ with associated isotropic vector $v_P=(n,p,q)^t$. Then the exhaustion function $f_P$ can be written as
\begin{align}\label{eq:exhaustion}
 f_{0}(X)&:= f_{P_0}(X)=\sqrt{-\frac{\sqrt{|D|}}{2}\QQC(X,X)}=y\\
 f_P(X)&=\frac{f_0(X)}{|\QQC(\sqrt{D}X,v_P)|}=\frac{y}{|q\bar{z}+\sqrt{D}p\bar{u}-n|}.
\end{align}

 These exhaustion functions are used to define the decomposition of $\D_0$ into sets $\D'(\I)$ for $\I \subset \Par$.  Let $\D(\I)\subset \D$ to be the set of $X \in \D$ such that $ f_P(X) \geq  f_Q(X)$ for every $P \in \I$ and $Q \in \Par \setminus \I$ and  $ f_P(X) =  f_P'(X)$ for every $P,P' \in \I$.  Define $\D'(\I)\subset \D'(\I)$ to be the subset where the inequality is strict, $ f_P(X) >  f_Q(X)$ for every $P \in \I$ and $Q \in \Par \setminus \I$. In other words, $\D'(\I)$ consists of the points that are higher with respect to $P\in \I$ than any other cusps.  The subset $\I$ is \emph{admissible} if $\D(\I) \neq \emptyset$ and \emph{strongly admissible} if $\D'(\I) \neq \emptyset$. 
\end{subsection}

\begin{subsection}{First contacts} \label{subsec:computing}
\begin{defn}
For a finite subset $\I \subset \Par$, the \emph{first contact for $\I$} is the subset of 
\[\{X \in \D\ |\ f_P(X)= f_Q(X)\text{ for every $P,Q \in \Par$} \}
\]
where the exhaustion functions $f_P$ attain their maximum.
\end{defn}

\begin{prop}\label{prop:first}
Let $P\neq P_0$ be a rational parabolic subgroup and let $v_P=(n,p,q)^t$ be the associated isotropic vector.  Then the first contact for $\{P_0,P\}$ is the point $X=(z,u)$, where
\[z=\frac{n}{q}+\frac{i}{|q|}\quad \text{and} \quad u=\frac{p}{q}.\]
In particular, \[f_0(X)^2=\frac{1}{|q|} =\frac{1}{\sqrt{|D|}|\QQC(v_0,v_P)|} .\]  More generally, the first contact $\hat{X}$ for rational parabolic subgroups $\{Q,R\}$ satisfies
\[f_Q(\hat{X})^2 = \frac{1}{\sqrt{|D|}|\QQC(v_Q,v_R)|}.\]
\end{prop}
\begin{proof}
To find the first contact for $\{P_0,P\}$, we need to find the point $X=(z,u)$ such that $f_0(X)$ is maximal on the set 
\[|\bar{z}+\sqrt{D}p\bar{u}/q-n/q|=1/|q|.\]  
Using \eqref{eq:coordinates} to express $z$ in its real and imaginary parts,
\[|r-iy^2-\sqrt{D}|u|^2/2+\sqrt{D}\bar{u}p/q-n/q|=1/|q|,\] 
we see that we can pick $r$ such that the real part of the expression is $0$.  Thus it suffices to consider points where 
\[|y^2+\sqrt{|D|}|u|^2/2-\Im(\sqrt{D}\bar{u}p/q)+\Im(n/q)|=1/|q|.\]
Since $v_P$ is isotropic, the left side of the equation is 
\[|y^2+\sqrt{|D|}|p|^2/(2|q|^2)+\sqrt{|D|}|u|^2/2-\Im(\sqrt{D}\bar{u}p/q)|=1/|q|,\]
and hence the first contact point occurs at the minimum of 
\[\sqrt{|D|}|u|^2/2-\Im(\sqrt{D}\bar{u}p/q)=\sqrt{|D|}(|u|^2/2-\Re(\bar{u}p/q))\]
One can compute that this implies that the first contact occurs when $u=p/q$.  Plugging back in, we get that $y=1/\sqrt{|q|}$ and $r=\Re(n/q)$, or equivalently, \[z=\frac{n}{q}+\frac{i}{|q|} \quad \text{as desired}.\] 

For the second statement, suppose $Q$ and $R$ are distinct rational parabolic subgroups. Let $\hat{X}$ denote the first contact for $\{Q,R\}$.  There exists $\gamma \in \Gb$ such that $\lsp{\gamma}Q = P_0$, $\lsp{\gamma}R = P$ for some $P$.  Note that $\gamma \cdot \hat{X}$ is the first contact for $\{P_0,P\}$.  By the $\Gb$-invariance of the exhaustion functions, it follows that 
\[f_Q(\hat{X})^2 =f_{P_0}(\gamma \cdot \hat{X})^2=\frac{1}{\sqrt{|D|}|\QQC(v_{P_0},v_P)|}.\]
Since $\Gb$ preserves $\QQC$, $|\QQC(v_{P_0},v_P)| = |\QQC(v_Q,v_R)|$ and the result follows.
\end{proof}
\end{subsection}
\begin{subsection}{Candidates}
In practice, once a lower bound $\mu$ is computed for the exhaustion functions (one exists because of reduction theory), one can compute a finite list subsets of rational parabolic subgroups from which a complete set of $\Gb$-representative of strongly admissible sets can be chosen.  

First, compute all the possible admissible sets of order 2 as follows.  Since $h(\Gb)=1$, every admissible set is $\Gb$-conjugate to an admissible set containing $P_0$.  Let $\I=\{P_0,P\}$ be such an admissible set.  Then by Proposition~\ref{prop:first}, if $P$ has associated isotropic vector $v_P=(n,p,q)^t$ then $|q| \leq 1/\mu^2$.  Since $q\in \OO$, this is a finite list of possibilities for $q$.  Furthermore, using $\Gb_{P_0}= \Gb \cap P_0$, we can move a point $X=(z,u)$ to a point $X'=(z',u')$, where $u'$ lies in a fundamental domain for $\CC$ modulo translation by $\OO$ and $f_0(X)=f_0(X')$.  In particular, we can arrange that the $u$-component of the first contact of $\I$ lies in a fundamental domain for $\CC$ modulo translation by $\OO$.  Since the first contact of $\{P_0,P\}$ has $u=p/q$, it follows that there are finitely many choices for $p$.  Finally, since $v_P$ is isotropic and there are only finitely many choices for $p$ and $q$, \eqref{eq:isotropic} implies there are finitely many choices for $n$.   

Once we have a finite list of candidates for $\Gb$-representatives of admissible sets of order two, we check if it is possible to find candidates of order greater than two that contain these candidates.  We repeatedly use Proposition~\ref{prop:first}.  In particular, if $\I$ is admissible, then  
\[|\QQC(v_Q,v_R)| \leq \frac{1}{\sqrt{|D|}\mu^2}\quad \text{for every $Q,R\in \I$.}\]  Given such an $\I$, we can check if it is admissible by seeing if $\D(\I) \neq \emptyset$.

We see that the isotropic vectors provide a convenient way to study admissible sets.  If the value of $|\QQC(v_P,v_Q)|$ is small  enough, this invariant is enough to classify the $\Gb$-conjugacy class of $\{P,Q\}$. 
\begin{prop}
If $\I=\{P,Q\}$ is a subset of rational parabolic subgroups such that $|\QQC(v_P,v_Q)|^2=1/|D|$, then $S$ is $\Gb$-equivalent to $\{P_0,P_w\}$, where $P_w$ is the rational parabolic subgroup that stabilizes the line through $(0,0,1)^t$. 
\end{prop}
\begin{proof}
Since the class number is $1$, every rational parabolic subgroup is $\Gb$-conjugate to $P_0$.  Without loss of generality, we can assume that $v_P=(1,0,0)^t$.  Then since $|\QQC(v_P,v_Q)|^2=1/|D|$, $v_Q$ is of the form $v=(n,p,1)$ for some $n,p \in \OO$.  We can multiply by any matrix in $\Gb_{N_0}=\Gb \cap N_0$ without changing $(1,0,0)^t$.  In particular,  
\[\mat{1&-\bar{p}\sqrt{D}&|p|^2\omega\\0&1&-p\\0&0&1}\vect{n\\p\\1} = \vect{n'\\0\\1}.\]  Note that since $(n',0,1)^t$ is isotropic, \eqref{eq:isotropic} implies that $n'\in \RR$. It follows that  $n' \in \ZZ$ and 
\[\mat{1&0&-n'\\0&1&0\\0&0&1} \in \Gb_{N_0}.\]  Thus we can arrange that $\{u,v\}$ is a $\Gb$-translate of a set of the form \[\{(1,0,0)^t,(0,0,1)^t\}.\]
\end{proof}
\end{subsection}
\end{section}

\begin{section}{Cohomology}\label{sec:cohomology}
We follow \cite{AGM} closely and refer to that paper for the details. Many of the complications arising from orientation issues disappear in our case since we subdivide to get a cell complex such that the stabilizer of each cell fixes the cell pointwise.  For the remainder of the paper, let $\G$ be a congruence group of level $N$ ($\G(N),\G_1(N),$ or $\G_0(N)$) and let $\Gb=\SU(2,1;\OO)$.  

\begin{subsection}{Cell complex}
The decomposition of the spine has a refinement $W$ into a regular cell complex such that the stabilizer of each cell fixes the cell pointwise.  We use $W$ to compute the integral cohomology of $\G \backslash \D$.

\begin{defn}
Let $\sigma$ be a cell of $W$.  Then the \emph{type of $\sigma$} is the $\Gb$-conjugacy class of $\sigma$.   
\end{defn}
\end{subsection}

\begin{subsection}{Orientation}
Let $\phi$ be an oriented cell of $W$ and let $\gamma \in \Gb$.  Since $\Gb$ acts by diffeomorphisms, $\gamma$ takes the orientation of $\phi$ to some orientation of $\gamma \phi$.  Let $\gamma_*\phi$ denote the cell $\gamma \phi$ with this choice of orientation.  

Let $W_T$ denote the $\Gb$-orbit of cells of type $T$.  For each type $T$, fix a representative cell $\phi_T \in W_T$.  Fix orientations on the standard cells $\phi_T$.  We use these orientations to fix orientations of all the cells of $W$.  In particular, if $\phi$ is a cell of type $T$ and $\phi=\gamma \phi_T$, then we give $\phi$ the orientation $\gamma_*\phi$.  Note that because the stabilizer of a cell fixes the cell pointwise, this gives a well-defined orientation to each cell of $W$.    
\end{subsection}
\begin{subsection}{Incidence function}
We now define the incidence function on $W$.  If $\psi$ is a facet of $\phi$, then the orientation of $\phi$ induces an orientation of $\psi$.  This orientation of $\psi$ may or may not agree with the orientation determined above.  Let $[\phi,\psi]$ be $\pm1$ depending on whether this orientation on $\phi$ does or does not agree with the orientation of $\psi$. This is precisely the $(\phi,\psi)$-component of the cellular boundary operator on $W$.  The incidence function $[\cdot,\cdot]$ is $\Gb$-invariant in the sense that $[\gamma_* \phi, \gamma_* \psi] =[\phi,\psi]$ for all $\gamma \in \Gb$. 

For every subgroup $\G \subset \Gb$, the quotient $\G \backslash W$ has the structure of a regular cell complex with identifications.  Since the incidence function defined above is $\Gb$-invariant, it is in particular $\G$-invariant.  Thus we can use $[\cdot,\cdot]$ to compute the cellular cohomology of $\G\backslash W$.  If $\phi$ is a cell of $W$, we will write $\overline{\phi}$ to denote its $\G$-conjugacy class, or equivalently the corresponding cell in $\G \backslash W$.
\end{subsection}
\begin{subsection}{Cells as orbits}  Let $X$ denote the homogeneous space $\G \backslash \Gb$.  The following proposition allows us to identify cells of the quotient $\G \backslash W$ with right $\Gb_T$-orbits in $X$.
\begin{prop}[{\cite[Proposition~3.3]{AGM}}]\label{prop:correspondence}
There is a one-to-one correspondence between 
\begin{enumerate}
\item Cells in $W_T$ and cosets $\Gb/\Gb_T$.
\item Cells of the quotient $\G\backslash W_T$ and $\G$-orbits of $W_T$.
\item Cells of the quotient $\G\backslash W_T$ and double cosets in $\G \backslash \Gb / \Gb_T$.
\end{enumerate}
\end{prop}
\begin{proof}
The first statement is a standard fact about stabilizers.  Since the stabilizer of a cell fixes the cell pointwise, the cell structure on $W$ descends to give a cell structure on $\G\backslash W$, and the second statement follows.  Combining the two statements gives that the cells of $\G\backslash W_T$ are parametrized by $\G \backslash\Gb/\Gb_T$.
\end{proof}

In light of this proposition, we will now think of a cell of $\G \backslash W$ as a $\G_T$-orbits of $X$ or a double coset in $\G \backslash \Gb/\Gb_{T}$.  

The faces of $\phi_T$ can be expressed as translates of the various representative cells $\phi_{T'}$.  In particular, there exist finite subsets $\BC(T,T') \subset \Gb$ such that 
\[\FF_{\phi_T} = \bigcup_{T'}\bigcup_{\alpha \in B(T,T')} \alpha \phi_{T'}.\]  
Note that the $\alpha$'s are not uniquely determined, but the cosets $\alpha \Gb_{T'} \in \Gb/\Gb_{T'}$ are determined for each pair $(T,T')$.  

The cell $\overline{\phi_T}$ corresponds to the double coset $\G e \Gb_T$, and the face  $\overline{\alpha \phi_{T'}}$ corresponds to the double coset $\G \alpha \Gb_{T'}$.  In terms of orbits, if $O$ is the $\Gb_{T}$-orbit corresponding to $\overline{\phi_T}$, then the orbit corresponding to $\overline{\alpha \phi_{T'}}$ is the $\Gb_{T'}$ orbit of $X$ that contains $O \cdot \alpha$.  By translation, we can understand the boundary faces of an arbitrary cell $\overline{\phi}$ of $\G \backslash W$ in terms of the $\alpha$ defined above.  

\begin{prop}[{\cite[Proposition~3.20]{AGM}}]\label{prop:orbit}
Let $\phi$ be a cell of type $T$ and let $O$ denote the $\Gamma_T$-orbit of $X$ corresponding to $\overline{\phi}$.  Then the faces of $\overline{\phi}$ of type $T'$ correspond exactly to the $\Gb_{T'}$-orbits of $X$ containing $O\cdot \alpha$ as $\alpha$ ranges over $\BC(T,T')$.
\end{prop}

From Propositions~\ref{prop:correspondence} and \ref{prop:orbit}, it becomes clear that to compute cohomology in this way, we need cellular structure of $W$, the stabilizers $\Gb_T$, the elements $\alpha \in \BC(T,T')$, and the incidence numbers $[\phi_T,\psi]$.  This is carried out for $k=\QQ(i)$ and $k=\QQ(\sqrt{-3})$ in the following sections.  
\end{subsection}
\end{section}

\begin{section}{Examples}\label{sec:examples}
\begin{subsection}{$k=\QQ(i)$}
In \cite{Yaspicard}, we compute the space $W$ for $k=\QQ(i)$.  The cells of $W$ fall into twenty-four equivalence classes modulo $\Gb$ consisting of two $3$-cells, seven $2$-cells, nine $1$-cells, and six $0$-cells.  Representatives of the $3$-cells and their boundary faces are shown in Figure~\ref{fig:3cells}.  We refer to \cite{Yaspicard} for the cell data.  The results of the cohomology computation are tabulated in Section~\ref{sec:numerical}.

\begin{figure}
\[\begin{array}{c @{\hspace{0.1in}} c}
{\xy/r.7in/: ="A", +(.2,1.1)="B","A",
{\xypolygon6"C"{~:{(1,-.2):(0,.5)::}~={30}{\bullet}}},
"B",{\xypolygon5"D"{~:{(.8,-.2):(0,.5)::}{\bullet}}},
"D2";"D4"**@{-},"D2";"D5"**@{-},
"D5";"C6"**@{-},
"D5";"C5"**@{-},
"D4";"C5"**@{-},
"D4";"C4"**@{-},
"D3";"C3"**@{-},
"D2";"C2"**@{-},
"D1";"C1"**@{-}
\endxy}
&
\xy/r.7in/: ="A", +(.2,1.1)="B","A",
{\xypolygon4{~:{(.8,-.2):(0,.7)::}~<>{;"B"**@{-}}{\bullet}}},
"B",\drop{\bullet}\endxy
\end{array}\]
\caption{The $3$-cells for $\SU(2,1;\ZZ[i])$.}
\label{fig:3cells}
\end{figure}
\end{subsection}

\begin{subsection}{$k=\QQ(\sqrt{-3}$)}\label{sec:eisenstein}
In this section we apply the previous results in the case $k=\QQ(\sqrt{-3})$.  Let $\zeta=(1+\sqrt{-3})/2$.  Let     
\begin{gather}
w=\begin{pmatrix} 0 & 0 & -1 \\ 0 & 1 & 0 \\ 1 & 0 & 0 \end{pmatrix} ,\quad \tau= \begin{pmatrix} 1 & 0 & 1 \\ 0 & 1 & 0 \\ 0 & 0 & 1 \end{pmatrix},\\
\gamma_1=\mat{\zeta^5&0&\zeta^2\\-1&\zeta^2&1\\\zeta^4&\sqrt{-3}\zeta&1}
, \quad \gamma_2=\mat{-1&\sqrt{-3}&\zeta^5\\0&1&-1\\-1&\sqrt{-3}&\zeta^4},\\
\epsilon=\mat{\zeta&0&0\\0&\zeta^{-2}&0\\0&0&\zeta},\quad \text{and}\quad \sigma = \mat{1&\sqrt{-3}&\zeta\\0&1&1\\0&0&1}.\end{gather}

A weak lower bound for the exhaustion functions is found by examining the values of various exhaustion functions on a set $\Omega \subset \D$ whose $\Gb$ translates cover $\D$.  Then Section~\ref{subsec:computing} is then used to compute the possible admissible sets.  Each was either rejected or verified using MAPLE.  The results are given below.

\begin{prop}\label{prop:types}
Every cell of $\D_0$ is $\G$-conjugate to exactly one of the following:
\begin{enumerate}
\item $\I^2=\left\{\vect{1\\0\\0},\vect{0\\0\\1}\right\}$
\item $\I^3_1=\left\{\vect{1\\0\\0},\vect{0\\0\\1},\vect{1\\0\\1}\right\}$
\item $\I^3_2=\left\{\vect{1\\0\\0},\vect{0\\0\\1},\vect{\zeta\\1\\1}\right\}$
\item $\I^4=\left\{\vect{1\\0\\0},\vect{0\\0\\1},\vect{1\\0\\1},\vect{\zeta\\1\\1}\right\}$
\item $\I^8=\left\{\vect{1\\0\\0},\vect{0\\0\\1},\vect{1\\0\\1},\vect{\zeta\\1\\1},\vect{\zeta^2\\1\\1},\vect{\zeta\\\zeta^5\\1},\vect{\sqrt{-3}\zeta^5\\\zeta^5\\1},\vect{\zeta^2\\\zeta\\1}\right\}$
\end{enumerate}
\end{prop}

\begin{table}
\caption{Incidence types}\label{tab:incidence}
\begin{tabular}{|c|c|cc|c|c|}
\hline
        &$\I^2$&$\I^3_1$&$\I^3_2$&$\I^4$ &$\I^8$\\\hline
$\I^2_1$&$\ast$&3       &3       &6     &28     \\\hline
$\I^3_1$&2     &$\ast$  &$\ast$  &1     &8      \\
$\I^3_2$&12    &$\ast$  &$\ast$  &3     &48     \\\hline
$\I^4$  &24    &6       &3       &$\ast$&16     \\\hline
$\I^8$  &12    &6       &3       &2     &$\ast$ \\
\hline
\end{tabular}
\end{table}
The incidence table is given in Table~\ref{tab:incidence}, where the entry below the diagonal means that each column cell has that many row cells in its boundary, and the entry above the diagonal means the column cell appears in the boundary of this many row cells.  The entries below the diagonal can be read off from Figure~\ref{fig:topcell}.  The entries above the diagonal can be easily computed from Proposition~\ref{prop:types}, since the $\Gamma$-conjugacy class of a strongly admissible can be distinguished by the pairwise $\mathcal{Q}$-inner products of its associated isotropic vectors, except to distinguish $\I^3_1$ and $\I^3_2$ we must also compute the dimension of the span of their isotropic vectors.  

Since the set of isotropic vectors associated with each piece of $\D_0$ is known, computing the stabilizers is a finite computation.  Let $\I$ be a strongly admissible set and let $\J$ be the associated set of isotropic vectors.  If $\gamma \in \Stab(\D(\I))$, then
\[\gamma \cdot \D(\I) = \D(\lsp{\gamma}\I) = \D(\I).\]
It follows that $\gamma$ must take each vector of $\J$ to another vector of $\J$ up to scaling by $\OO^*$. 

\begin{prop}
The stabilizer of $\I^2$ is isomorphic to $\ZZ/12\ZZ$ and is generated by $\epsilon w$.  The stabilizer of $\I^3_1$ is isomorphic to $\ZZ/3\ZZ \oplus \ZZ/6\ZZ$ and is generated by $\tau w$ and $\epsilon$.  
The stabilizer of $\I^3_2$ is isomorphic to $\ZZ/3\ZZ$ and is generated by $\epsilon^2$.  The stabilizer of $\I^4$ is isomorphic to $\ZZ/3\ZZ$ and is generated by $\epsilon^2$.  The stabilizer $\I^8$ is the order $24$ group, Magma small group number $11$ and is generated by $\gamma_1$ and $\gamma_2$.

\end{prop}
\begin{figure}
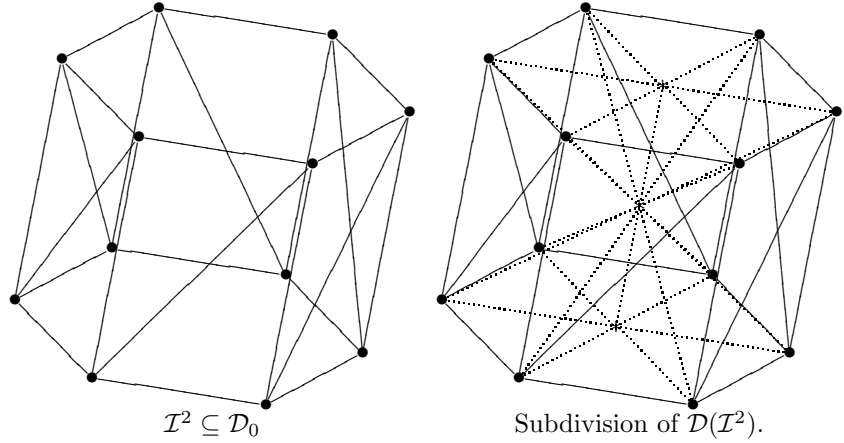

\[
\begin{array}{c @{\hspace{0.1in}} c}
{\xy/r.7in/: ="A", +(.35,1.8)="B","A",
{\xypolygon6"C"{~:{(1.3,-.2):(0,.43)::}{\bullet}}},
"B",{\xypolygon6"D"{~:{(1.3,-.2):(0,.43)::}{\bullet}}},
"D5";"C5"**@{-},
"D6";"C5"**@{-},
"D6";"C6"**@{-},
"D1";"C6"**@{-},
"D1";"C1"**@{-},
"D2";"C1"**@{-},
"D2";"C2"**@{-},
"D3";"C2"**@{-},
"D3";"C3"**@{-},
"D4";"C3"**@{-},
"D4";"C4"**@{-},
"D5";"C4"**@{-},
\endxy} & 
{\xy/r.7in/: ="A", +(.35,1.8)="B",-(.175,0.9)="E","A",
{\xypolygon6"C"{~:{(1.3,-.2):(0,.43)::}~<{.}{\bullet}}},
"B",{\xypolygon6"D"{~:{(1.3,-.2):(0,.43)::}~<{.}{\bullet}}},
"D5";"C5"**@{-},
"D6";"C5"**@{-},
"D6";"C6"**@{-},
"D1";"C6"**@{-},
"D1";"C1"**@{-},
"D2";"C1"**@{-},
"D2";"C2"**@{-},
"D3";"C2"**@{-},
"D3";"C3"**@{-},
"D4";"C3"**@{-},
"D4";"C4"**@{-},
"D5";"C4"**@{-},
"D0";"C0"**@{.},
"E";"C1"**@{.},
"E";"C2"**@{.},
"E";"C3"**@{.},
"E";"C4"**@{.},
"E";"C5"**@{.},
"E";"C6"**@{.},
"E";"D1"**@{.},
"E";"D2"**@{.},
"E";"D3"**@{.},
"E";"D4"**@{.},
"E";"D5"**@{.},
"E";"D6"**@{.},
"E"*{\ast},
"D0"*{\ast},
"C0"*{\ast}
\endxy}\\
\I^2 \subseteq \D_0&\text{Subdivision of $\D(\I^2)$.}
\end{array}
\]
\caption{The $3$-cell $\D(\I^2)$ and subdivision for $\SU(2,1;\ZZ[\zeta])$.}\label{fig:topcell}
\end{figure}

Since the action of $\Gb$ on $\D_0$ is understood, we can subdivide $\D(\I^2)$ and its boundary into cells in a $\Gb$-equivariant way as shown in Figure~\ref{fig:topcell} to yield a cell complex $W$ with the property that the stabilizer of each cell fixes the cell pointwise.  The cells of $W$ fall into sixteen equivalence classes modulo $\Gb$ consisting of one $3$-cells, seven $2$-cells, six $1$-cells, and three $0$-cells.  Representatives are chosen and stabilizers are recomputed.  The results are given in Table~\ref{tab:stabilizer}. 

\begin{table}
\caption{Stabilizers of representative cells}\label{tab:stabilizer}
\begin{tabular}{|c|c|c|c|}
\hline
Cell &Dimension&Stabilizer& Generators\\\hline\hline
$X$&$3$&$\ZZ/3\ZZ$&$\epsilon^2$\\\hline
$A$&$2$&$\ZZ/3\ZZ$&$\epsilon^2$\\
$B$&$2$&$\ZZ/3\ZZ$&$\epsilon^2$\\
$C$&$2$&$\ZZ/3\ZZ$&$\epsilon^2$\\
$E$&$2$&$\ZZ/3\ZZ$&$\epsilon^2$\\
$F$&$2$&$\ZZ/3\ZZ$&$\epsilon^2$\\
$G$&$2$&$\ZZ/3\ZZ$&$\epsilon^2$\\\hline
$a$&$1$&$\ZZ/3\ZZ$&$\epsilon^2$\\
$b$&$1$&$\ZZ/3\ZZ$&$\epsilon^2$\\
$c$&$1$&$\ZZ/3\ZZ$&$\epsilon^2$\\
$d$&$1$&$\ZZ/3\ZZ$&$\epsilon^2$\\
$e$&$1$&$\ZZ/3\ZZ$&$\epsilon^2$\\
$f$&$1$&$\ZZ/6\ZZ$&$\epsilon$\\\hline
$m$&$0$&$\GC(24,11)$\footnotemark&$\gamma_1, \gamma_2$\\
$n$&$0$&$\ZZ/12\ZZ$&$\epsilon w$\\
$o$&$0$&$\ZZ/3\ZZ \oplus \ZZ/6\ZZ$&$\epsilon^2,\tau \epsilon w$\\
\hline
\end{tabular}
\end{table}
\footnotetext{This is the order $24$ with Magma small group library number $11$.}
\end{subsection}
\end{section}

\begin{section}{Numerical results}\label{sec:numerical}
A computer was used to compute the $\Gb_T$-orbits for the various congruence subgroups $\G$ and to create the coboundary matrix.  In this process, large sparse matrices were produced, up to $97,542 \times 89,478$ for $\G_0(-2+13\sqrt{-3}) \subset \SU(2,1;\ZZ[\zeta])$.  Then MAGMA \cite{magma} was used to compute the elementary divisors.  The computations were only carried out in a few cases for $\G_1(N)$ and $\G(N)$ because the index grows so rapidly.  In fact, Holzapfel \cite{Hol} has an explicit formula for the index $[\G(N):\Gb]$ for $N \in \ZZ$, which shows that the index grows like $N^8$.  Tables~\ref{tab:coh0inert}-\ref{tab:cohprinc} summarize the numerical results for $k=\QQ(i)$ and Table~\ref{tab:eis} summarizes the results for $k=\QQ(\sqrt{-3})$.  We remark that although the cohomology in degree $1$ is predominately $0$ in the range that was computed, by a result of Blasius and Rogawski \cite{BlRo}, we expect infinitely many congruence subgroups to yield non-trivial cohomology.  Their result, however, is for principal congruence subgroups, not those of the form $\G_0(N)$, and we already see some non-trivial cohomology in $H^1$ for the former.

\begin{table}
\caption{The integral cohomology for $\G_0(N) \subset \SU(2,1;\ZZ[i])$}\label{tab:coh0inert}
\begin{tabular}{|r||l|l|l|l|}
\hline
$N$  & $H^1$ & $H^2$&$H^3$\\
\hline
$2$&$0$&$\ZZ$&$\ZZ^2$\\ 
$3$&$0$&$\ZZ^2$&$\ZZ \oplus \ZZ/2\ZZ$\\ 
$4$&$0$&$\ZZ^2\oplus \ZZ/2\ZZ$&$\ZZ^6$\\ 
$5$&$0$&$\ZZ^7$&$\ZZ^5$\\ 
$6$&$0$&$\ZZ^4 \oplus \ZZ/2\ZZ$&$\ZZ^5\oplus \ZZ/2\ZZ$\\ 
$7$&$0$&$\ZZ^{10}\oplus \ZZ/4\ZZ$&$\ZZ\oplus \ZZ/12\ZZ$\\
$8$&$0$&$\ZZ^{11} \oplus (\ZZ/2\ZZ)^3 \oplus \ZZ/4\ZZ $&$ \ZZ^{12} $\\
$9$&$0$&$\ZZ^{25} \oplus (\ZZ/3\ZZ)^2 $&$\ZZ^5\oplus (\ZZ/2\ZZ)^2  $\\
$10$&$0$&$\ZZ^{34}$&$\ZZ^{17}$\\
$11$&$0$&$\ZZ^{39}\oplus \ZZ/15\ZZ$&$\ZZ \oplus \ZZ/30\ZZ$\\
$12$&$0$&$\ZZ^{36}\oplus \ZZ/2\ZZ \oplus \ZZ/8\ZZ$&$\ZZ^{13} \oplus \ZZ/13\ZZ$\\
$13$&$0$&$\ZZ^{79}$&$\ZZ^{5} \oplus (\ZZ/3\ZZ)^2$\\
$14$&$0$&$\ZZ^{58}\oplus \ZZ/12\ZZ $&$\ZZ^{5} \oplus \ZZ/12\ZZ$\\
$15$&$0$&$\ZZ^{148}\oplus \ZZ/2\ZZ $&$\ZZ^{11} \oplus (\ZZ/2\ZZ)^3$\\
$16$&$0$&$\ZZ^{94}\oplus (\ZZ/2\ZZ)^3\oplus (\ZZ/4\ZZ)^4 \oplus (\ZZ/8\ZZ)^2 $&$\ZZ^{23} \oplus \ZZ/2\ZZ$\\
$17$&$0$&$\ZZ^{166}\oplus (\ZZ/2\ZZ)^2$&$\ZZ^{5} \oplus (\ZZ/2\ZZ)^2 \oplus (\ZZ/4\ZZ)^2$\\
$18$&$0$&$\ZZ^{142}\oplus \ZZ/3\ZZ \oplus Z/6\ZZ$&$\ZZ^{17} \oplus (\ZZ/2\ZZ)^3 $\\
$19$&$0$&$\ZZ^{211} \oplus \ZZ /15 \ZZ$ &$\ZZ \oplus \ZZ/90\ZZ$\\ 
$20$&$0$&$\ZZ^{238}\oplus (\ZZ/2\ZZ)^{5}$ &$\ZZ^{41} \oplus (\ZZ/2\ZZ)^2$\\ 
$21$&$0$&$\ZZ^{294}\oplus (\ZZ/2\ZZ)^{2}\oplus \ZZ/16\ZZ$ &$\ZZ^{3} \oplus (\ZZ/2\ZZ)^2\oplus \ZZ/48\ZZ$\\ 
$22$&$0$&$\ZZ^{238}\oplus \ZZ/30\ZZ$ &$\ZZ^{5} \oplus \ZZ/30\ZZ$\\ 
$23$&$0$&$\ZZ^{372}\oplus \ZZ/132\ZZ$&$\ZZ \oplus \ZZ/132\ZZ$\\
$24$&$0$&$\ZZ^{312}\oplus (\ZZ/2\ZZ)^5\oplus \ZZ/4\ZZ \oplus \ZZ/8\ZZ$&$\ZZ^{25} \oplus (\ZZ/2\ZZ)^2\oplus \ZZ/8\ZZ$\\
\hline
$1+i$&$0$&$\ZZ$&$\ZZ$\\\hline 
$2+i$&$0$&$\ZZ^3$&$\ZZ^2$\\ 
$3+2i$&$0$&$\ZZ^7$&$\ZZ^2\oplus \ZZ/3\ZZ$\\ 
$4+i$&$0$&$\ZZ^9\oplus \ZZ/2\ZZ$&$\ZZ^2\oplus \ZZ/2\ZZ \oplus \ZZ/4\ZZ$\\
$5+2i$&$0$&$\ZZ^{21}\oplus \ZZ/7\ZZ$&$\ZZ^2\oplus \ZZ/7\ZZ$\\
$6+i$&$0$&$\ZZ^{37}\oplus \ZZ/3\ZZ$&$\ZZ^2\oplus \ZZ/9\ZZ$\\
$5+4i$&$0$&$\ZZ^{45}\oplus \ZZ/5\ZZ$&$\ZZ^2\oplus \ZZ/2\ZZ \oplus \ZZ/10\ZZ$\\
$7+2i$&$0$&$\ZZ^{75}\oplus \ZZ/13\ZZ$&$\ZZ^2\oplus \ZZ/13\ZZ$\\
$6+5i$&$0$&$\ZZ^{103} \oplus \ZZ/5\ZZ$&$\ZZ^2\oplus \ZZ/15 \ZZ$\\
$8+3i$&$0$&$\ZZ^{151}\oplus \ZZ/3\ZZ$&$\ZZ^2\oplus \ZZ/2\ZZ \oplus \ZZ/18\ZZ$\\
$8+5i$&$0$&$\ZZ^{225}\oplus \ZZ/11\ZZ$&$\ZZ^2\oplus \ZZ/2\ZZ \oplus \ZZ/22\ZZ$\\
$9+4i$&$0$&$\ZZ^{271}\oplus \ZZ/4\ZZ$&$\ZZ^2\oplus \ZZ/2\ZZ \oplus \ZZ/24\ZZ$\\
$10+i$&$0$&$\ZZ^{291}\oplus \ZZ/25\ZZ$&$\ZZ^2\oplus \ZZ/25\ZZ$\\
$10+3i$&$0$&$\ZZ^{343}\oplus \ZZ/9\ZZ$&$\ZZ^2\oplus \ZZ/27\ZZ$\\
$8+7i$&$0$&$\ZZ^{369}\oplus \ZZ/14\ZZ$&$\ZZ^2\oplus \ZZ/2\ZZ \oplus \ZZ/28\ZZ$\\
\hline
\end{tabular}
\end{table}

\begin{table}
\caption{The integral cohomology of $\G_1(N) \subset \SU(2,1;\ZZ[i])$}\label{tab:coh1}
\begin{tabular}{|r|r||l|l|l|l|}
\hline
$N$ & $H^1$ & $H^2$&$H^3$\\
\hline
$2$&$0$&$\ZZ \oplus \ZZ/2\ZZ$&$\ZZ^2$\\
$3$&$0$&$\ZZ^{13}\oplus (\ZZ/2\ZZ)^2 \oplus \ZZ/4\ZZ$&$\ZZ^3$\\
$4$&$\ZZ^2$&$\ZZ^{24} \oplus (\ZZ/2\ZZ)^3\oplus (\ZZ/4\ZZ)^2\oplus \ZZ/8\ZZ$&$\ZZ^{11}$\\
$5$&$0$&$\ZZ^{115} \oplus (\ZZ/2\ZZ)^5\oplus (\ZZ/10\ZZ)^2\oplus \ZZ/20\ZZ$&$\ZZ^{23}$\\
$6$&$0$&$\ZZ^{102} \oplus (\ZZ/2\ZZ)^5\oplus \ZZ/6\ZZ\oplus (\ZZ/12\ZZ)^2$&$\ZZ^{19}$\\
$7$&$0$&$\ZZ^{538}\oplus (\ZZ/2\ZZ)^2\oplus (\ZZ/14\ZZ)^2 \oplus \ZZ/28\ZZ$&$\ZZ^{23}$\\
$8$&$\ZZ^6$&$\ZZ^{460} \oplus (\ZZ/2\ZZ)^5 \oplus (\ZZ /4\ZZ)^4\oplus \ZZ/8\ZZ \oplus \ZZ/16\ZZ$&$\ZZ^{71} \oplus \ZZ/2\ZZ$\\
\hline
$2+i$&$0$&$\ZZ^{10}\oplus \ZZ/2\ZZ$&$\ZZ^{3}$\\
$3+2i$&$0$&$\ZZ^{91}\oplus \ZZ/2\ZZ$&$\ZZ^{11}$\\
$4+i$&$0$&$\ZZ^{184}\oplus (\ZZ/2\ZZ)^2$&$\ZZ^{15}\oplus \ZZ/2\ZZ$\\
\hline
\end{tabular}
\end{table}

\begin{table}
\caption{The integral cohomology for $\G(N) \subset\SU(2,1;\ZZ[i])$}\label{tab:cohprinc}
\begin{tabular}{|c|r||l|l|l|l|}
\hline
$N$ & $H^1$ & $H^2$&$H^3$\\
\hline
$1+i$&$0$&$\ZZ$&$\ZZ^2$\\
$2$&$0$&$\ZZ^2\oplus (\ZZ/2\ZZ)^2$&$\ZZ^5$\\
$3$&$0$&$\ZZ^{243}\oplus (\ZZ/2\ZZ)^7\oplus (\ZZ/6\ZZ)^7 \oplus \ZZ/12\ZZ$&$\ZZ^{55}$\\
$4$&$\ZZ^6$&$\ZZ^{484}\oplus (\ZZ/2\ZZ)^7 \oplus (\ZZ/4\ZZ)^{12} \oplus (\ZZ/8\ZZ)^2$&$\ZZ^{95} \oplus \ZZ/2\ZZ$\\
\hline
\end{tabular}
\end{table}

\begin{table}
\caption{The integral cohomology for $\G_0(N) \subset \SU(2,1;\ZZ[\zeta])$}\label{tab:eis}
$
\begin{array}{|c|r||l|l|l|l|}
\hline
N & H^1 & H^2&H^3\\
\hline
5 & 0 & \ZZ^8 \oplus \ZZ/6\ZZ&\ZZ \oplus \ZZ/4\ZZ\\
11 & 0 & \ZZ^{56} \oplus (\ZZ/2\ZZ)^2 \oplus \ZZ/60\ZZ& \ZZ \oplus \ZZ/20\ZZ\\
17 & 0 & \ZZ^{208} \oplus \ZZ/3\ZZ \oplus \ZZ/24\ZZ& \ZZ \oplus \ZZ/48\ZZ\\
23 & 0 & \ZZ^{508} \oplus (\ZZ/2\ZZ)^2\oplus \ZZ/264\ZZ & \ZZ \oplus \ZZ/88\ZZ\\
\hline
2+\sqrt{-3}& 0 & \ZZ^7 & \ZZ^2\\ 
3+\sqrt{-3}& 0 & \ZZ^{15} & \ZZ^2 \oplus \ZZ/2\ZZ\\
-2+5\sqrt{-3}& 0&\ZZ^{21}\oplus \ZZ/3\ZZ&\ZZ^2 \oplus \ZZ/3\ZZ\\
5+\sqrt{-3}& 0&\ZZ^{47} \oplus \ZZ/5\ZZ&\ZZ^2 \oplus \ZZ/5\ZZ\\
-3+7\sqrt{-3}&0&\ZZ^{67} \oplus \ZZ/3\ZZ&\ZZ^2 \oplus \ZZ/6\ZZ\\
6+\sqrt{-3}&0&\ZZ^{85} \oplus \ZZ/7\ZZ & \ZZ^2 \oplus \ZZ/7\ZZ\\
-4+9\sqrt{-3}&0&\ZZ^{167} \oplus \ZZ/5\ZZ & \ZZ^2 \oplus \ZZ/10\ZZ\\
-2+9\sqrt{-3}&0&\ZZ^{197} \oplus \ZZ/11\ZZ & \ZZ^2 \oplus \ZZ/11\ZZ\\
8+\sqrt{-3} & 0 & \ZZ^{235} \oplus \ZZ/6\ZZ& \ZZ^2 \oplus \ZZ/12\ZZ\\
3+7\sqrt{-3}&0&\ZZ^{271} \oplus \ZZ/13 \ZZ & \ZZ^2 \oplus \ZZ/13\ZZ\\
-3+11\sqrt{-3}&0&\ZZ^{408} \oplus \ZZ/8\ZZ & \ZZ^2 \oplus \ZZ/16\ZZ\\
2+9\sqrt{-3}& 0 & \ZZ^{405} \oplus \ZZ/17\ZZ & \ZZ^2 \oplus \ZZ/17\ZZ\\
7+5\sqrt{-3}& 0 & \ZZ^{511} \oplus \ZZ/9\ZZ &\ZZ^2 \oplus \ZZ/18\ZZ\\
-6+13\sqrt{-3}&0&\ZZ^{687} \oplus \ZZ/21\ZZ & \ZZ^2 \oplus \ZZ/21\ZZ\\
\hline
\end{array}
$
\end{table}
\end{section}

\bibliography{../../references}    

\providecommand{\bysame}{\leavevmode\hbox to3em{\hrulefill}\thinspace}
\providecommand{\MR}{\relax\ifhmode\unskip\space\fi MR }
\providecommand{\MRhref}[2]{%
  \href{http://www.ams.org/mathscinet-getitem?mr=#1}{#2}
}
\providecommand{\href}[2]{#2}
\begin{thebibliography}{10}

\bibitem{A2}
A.~Ash, \emph{Small-dimensional classifying spaces for arithmetic subgroups of
  general linear groups}, Duke Math. J. \textbf{51} (1984), no.~2, 459--468.

\bibitem{AGM}
A.~Ash, P.~E. Gunnells, and M.~McConnell, \emph{Cohomology of congruence
  subgroups of {${\rm SL}\sb 4(\mathbb Z)$}}, J. Number Theory \textbf{94}
  (2002), no.~1, 181--212.

\bibitem{BlRo}
D.~Blasius and J.~Rogawski, \emph{Cohomology of congruence subgroups of {${\rm
  SU}(2,1)\sp p$} and {H}odge cycles on some special complex hyperbolic
  surfaces}, Regulators in analysis, geometry and number theory, Progr. Math.,
  vol. 171, Birkh\"auser Boston, Boston, MA, 2000, pp.~1--15.

\bibitem{magma}
W.~Bosma, J.~Cannon, and C.~Playoust, \emph{The {M}agma algebra system. {I}.
  {T}he user language}, J. Symbolic Comput. \textbf{24} (1997), no.~3-4,
  235--265, Computational algebra and number theory (London, 1993).

\bibitem{De}
L.~Demb{\'e}l{\'e}, \emph{Examples of automorphic forms on the unitary group
  {$U(3)$}}, Preprint, 2007.

\bibitem{FaPa}
E.~Falbel and J.~R. Parker, \emph{The geometry of the {E}isenstein-{P}icard
  modular group}, Duke Math. J. \textbf{131} (2006), no.~2, 249--289.

\bibitem{Fra}
J.~Franke, \emph{Harmonic analysis in weighted {$L\sb 2$}-spaces}, Ann. Sci.
  \'Ecole Norm. Sup. (4) \textbf{31} (1998), no.~2, 181--279.

\bibitem{Hol}
R.-P. Holzapfel, \emph{Zeta dimension formula for {P}icard modular cusp forms
  of neat natural congruence subgroups}, Abh. Math. Sem. Univ. Hamburg
  \textbf{68} (1998), 169--192.

\bibitem{Yasrank1}
D.~Yasaki, \emph{On the existence of spines for {$\mathbb Q$}-rank 1 groups},
  Selecta Math. (N.S.) \textbf{12} (2006), no.~3-4, 541--564.

\bibitem{Yaspicard}
\bysame, \emph{An explicit spine for the {P}icard modular group over the
  {G}aussian integers}, Journal of Number Theory (to appear), available online
  {\url{http://dx.doi.org/10.1016/j.jnt.2007.03.008}}.

\bibitem{Z}
T.~Zink, \emph{\"{U}ber die {A}nzahl der {S}pitzen einiger arithmetischer
  {U}ntergruppen unit\"arer {G}ruppen}, Math. Nachr. \textbf{89} (1979),
  315--320.

\end{thebibliography}
\end{document}